\documentclass[a4paper, 12pt] {article}

\usepackage{fender}  

\begin{document}

    \title{RCF\,1 \\
           Theories of $\mr{PR}$ Maps and Partial $\mr{PR}$ Maps \\
               $\mathbf{P\widehat{R}_A}$}
        
\footnotetext{
  This is part 1 of a cycle on \emph{Recursive Categorical Foundations.} \\ 
  There is a still more detailed version (pdf file) equally entitled
  \emph{Theories of PR Maps and Partial PR Maps}}

\footnotetext{
  Legend of LOGO: $\mathbf{P\widehat{R}_A} = \mathbf{\widehat{PR_A}}$
  symbolises the theory of \emph{partial} maps over categorical 
  Free-Variables theory $\PRa$ of Primitive Recursion, $\PR$, with 
  predicate Abstraction}

    \author{Michael Pfender\footnote{
              TU Berlin, Mathematik, pfender@math.tu-berlin.de}
           }    

    \date{May 2008\footnote{last revised \today}}                          

\maketitle

\textbf{Abstract:}
We give to the categorical theory $\PR$ of Primitive Recursion 
a logically simple, \emph{algebraic} presentation, via 
\emph{equations between maps,} plus one genuine \NAME{Horner} 
type schema, namely Freyd's uniqueness of the \emph{initialised iterated.} 
Free Variables are introduced -- formally -- as another 
names for \emph{projections.} Predicates $\chi: A \to 2$ admit 
\emph{interpretation} as (formal) Objects $\set{A\,|\,\chi}$ of a 
surrounding Theory $\PRa = \PR+(\abstr):$ schema $(\abstr)$ formalises
this \emph{predicate abstraction} into additional Objects. 
Categorical Theory $\hatPRa \bs{\sqsup} \PRa \bs{\sqsup} \PR$
then is the Theory of formally \emph{partial} $\mr{PR}$-maps, 
having Theory $\PRa$ embedded. This Theory $\hatPRa$ bears the 
structure of a (still) diagonal monoidal category. It is equivalent 
to ``the'' categorical theory of $\mu$-recursion (and of $\while$ loops), 
viewed as partial $\mr{PR}$ maps. So the present approach to partial 
maps sheds new light on \NAME{Church}'s \textsc{Thesis,} ``embedded'' 
into a Free-Variables, formally \emph{variable-free} (categorical) 
framework. 


\newpage

\section{Introduction} 

We develop here, from scratch, a formally 
\emph{variable-free, categorical} Theory of $\mu$-\emph{recursion,} 
without use of formal \emph{quantification:} 
This Theory is formalised on the basis of a theory $\hatPRa$ of 
\emph{partial} $\mr{PR}$ maps which in turn is introduced as a 
\emph{definitional,} conservative extension of Theory 
$\PRa = \PR+(\abstr),$ the latter obtained from \emph{fundamental} 
categorical Theory $\PR$ of Primitive Recursion: we formally 
interpret $\mr{PR}$ \emph{predicates} $\chi: A \to 2$ as additional, 
\emph{defined} Objects of emerging Theory $\PRa:$ schema of abstraction 
already ``hidden'' in \emph{fundamental} Theory $\PR$ of Primitive 
Recursion. The latter is given as \emph{Cartesian Hull} over
data and axioms of a Natural Numbers Object $\N,$ this in 
the sense of \NAME{Lawvere,} \NAME{Eilenberg} \& \NAME{Elgot,} and 
\NAME{Freyd.}

\smallskip
Central in present approach is the notion of a \emph{partial $\mr{PR}$ map:} 
Such a (formally) partial map $f: A \parto B$ is given, in categorical, 
variable-free terms, by:
\begin{itemize}
\item 
an \emph{enumeration source} Object $D_f$ in $\PRa$(\eg: $D_f = \N$),

\item a $\PRa$-map $d_f: D_f \to A,$ meaning for $\PRa$-\emph{enumeration} 
of \emph{defined arguments} of $f,$

\item
and a $\PRa$-\emph{rule}-map
$\widehat{f}: D_f \to B$ for $f,$ 

\item
these data with (intuitive) meaning: for $a \in A$ 
\emph{defined argument,} of form $a = d_f(\hat{a}) \in A$ 
for suitable $\hat{a} \in D_f: 
  \xymatrix{d_f(\hat{a}) = a \ar @{|->} [r]^<<<<<{f} 
            & \widehat{f} (\hat{a}) \in B}.$
\end{itemize}

\smallskip
We prove a \textbf{Structure Theorem} for 
\emph{partial-map extensions} $\widehat{\bfS} \bs{\sqsupset} \bfS,$ 
where $\bfS$ is a Cartesian $\mr{PR}$ theory with schema of predicate 
abstraction -- mostly $\bfS :\,= \PRa$ or one of its strengthenings -- 
which establishes these extensions (via \ul{embeddin}g) as 
\emph{diagonal monoidal} $\mr{PR}$ theories: Cartesian structure is lost in 
part, since the (still present) projection- and terminal-map-families 
do not preserve their character as \emph{natural transformations} 
in the extension. These \emph{partial-map extensions} turn out to be 
\emph{Closures:}  
$\hathatS \bs{\iso} \hatS:$ \emph{Partial} partial maps have a 
representation as just \emph{partial maps.}

Within this Free-Variables (formally: variable free) categorical 
framework $\hatPRa$ for \emph{partial} $\mr{PR}$ map theories, we discuss
(Free-Variables) \emph{category based} $\mu$-recursion as well as
\emph{content driven} loops such as $\while$ loops. 

This prepares in particular discussion of \emph{termination} for
suitable special such loops, namely those given as 
\emph{Complexity Controlled Iterations,}  
for which iteration \emph{step} decreases a \emph{complexity}
measure within a suitably given (constructive) \emph{ordinal} $O,$ 
``until'' minimum $0$ of $O$ is reached. Complexity Controlled 
Iteration is basic for the following second part of investigation
on \emph{Recursive Categorical Foundations,} RCF 2, entitled 
\emph{Evaluation and Consistency.}

Evaluation of $\mr{PR}$ map codes is there \emph{resolved} into such an
iteration with \emph{descending} complexity values. ``Hence'' we 
can ``hope'' this formally partial evaluation to always 
\emph{terminate,} by reaching complexity $0.$ In that case within 
ordinal $O$ taken the lexicographically ordered set of polynomials 
over $\N,$ in one indeterminate.



\section{Notions, Axioms,  Results for Theories $\PR$ and $\PRa$}

Fundamental Theory $\PR$ of \emph{Primitive Recursion} is the
\emph{minimal, ``initial''} \emph{Cartesian} Theory with 
(universal) Natural Numbers Object:

As Objects it has $\one,\ N,\ldots,A,\ B,\ldots,\ (A \times B),$ 
\ie all (binary bracketed) finite powers of Natural Numbers Object 
(``NNO'') $\N.$

It comes with associative \emph{composition} ``$\circ$'', 
functorial \emph{cylindrification} 
$(A \times g): A \times B \to A \times B'$ -- and from this with 
bifunctorial \emph{Cartesian product} ``$\times$'' -- and is 
generated over \emph{basic map constants} $0: \one \to \N$ 
\emph{(zero),} \emph{successor} $s: \N \times \N,$ as well as 
\emph{terminal map} $!: A \to \one,$ \emph{diagonal} 
$\Delta: A \to A \times A,$ and \emph{projections}
$\ell: A \times B \to A$ and $r: A \times B \to B.$

\smallskip
For given $f: C \to A$ and $g: C \to B,$ \emph{induced map} 
$(f,g): C \to A \times B$ into Cartesian Product is \textbf{defined} as
  $$(f,g) \defeq (f \times g) \circ \Delta_C: 
                      C \to C^2 \bydefeq C \times C \to A \times B,$$
with $(f \times g): C^2 \to A \times B$ defined by either -- equal(!) -- 
sequence of cylindrified above, see diamond sub-\textsc{diagram} in

\begin{minipage} {\textwidth}
$$
\xymatrix{
A 
\ar [rr]^{f}
& & A'
\\
& A' \times B
  \ar [ru]^{\ell}
  \ar [rd]^{A' \times g}
  &
\\
A \times B 
\ar [uu]^{\ell}
\ar [ru]^{f \times B}
\ar [rr]^{\bydefeq f \times g}
\ar [rd]_{A \times g}
\ar [dd]_r
& & A' \times B'
    \ar [uu]_{\ell}
    \ar [dd]^{r}
\\
& A \times B'
  \ar [ru]_{f \times B'}
  \ar [rd]_r
  &
\\
B 
\ar [rr]_g
& & B'
}
$$
\begin{center} Binary Cartesian Map Product \textsc{diagram} \end{center}
\end{minipage}

\bigskip



In present context we take \NAME{Godement}'s equations for Cartesian (!)
Product, as \textbf{axioms:}
$$
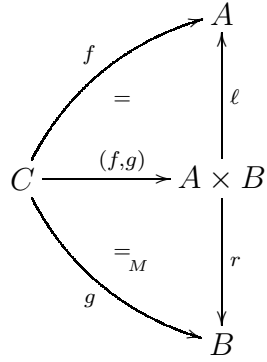
\xymatrix @+2em{
& A
\\
C
\ar @/^1pc/[ru]^{f}
\ar @{} [ru]| {\overset{\quad} {=}} 
\ar [r]^{(f,g)}
\ar @{} [rd]| {\underset{\quad{M}} {=}} 
\ar @/_1pc/ [rd]_{g}
& A \times B
  \ar [u]_{\ell}
  \ar [d]^{r}
\\
& B
}
$$

\textbf{Uniqueness} of the induced 
$(f,g) \bydefeq (f \times g)\,\circ\,\Delta_C: C \to C^2 \to A \times B$
is \emph{forced} -- \textbf{axiomatically} -- by the following 
\NAME{Fourman}'s equation (equational \emph{schema}):  
\inference{ (\mathrm{Four}!) }
{ $h: C \to A \times B$ map }
{ $(\ell_{A,B} \circ h,\ r_{A,B} \circ h) = h: C \to A \times B.$
}

\bigskip
The NNO $\xymatrix{
        \one
          \ar[r]^{0}
          & \N
            \ar[r]^{s}
            & \N}$      
is given -- by \textbf{axiom} -- the following property, expressed as 
\emph{iteration schema:}

\begin{minipage} {\textwidth}
$$
\xymatrix{
& A \times \N 
  \ar [r]^{A \times s}
  \ar [dd]^{f^\S}
  \ar @{} [ddr] |{=}
  & A \times \N 
      \ar [dd]^{f^\S}
      & 
\\
A
\ar [ur]^{(\id,0\,!)}
\ar [dr]_{\id}
\ar @{} [r] |{=}
& & 
  & (\mathrm{it})
\\
& A
  \ar[r]^f
  & A
    &
} 
$$
\begin{center} Basic Iteration \textsc{diagram} \end{center} 
\end{minipage}

\bigskip


\textbf{As Uniqueness schema} we need the following 
one -- of \NAME{Freyd,} formally stronger than 
\emph{uniqueness of the iterated} $f^\S: A \times \N \to A$ 
satisfying the equations of Basic Iteration 
\textsc{diagram} above -- namely 
\textbf{Uniqueness of initialised iterated:} 
\inference { (\FR!) }
{ $f: A \to B,\ g: B \to B,$ $h: A \times \N \to B,$ all in $\PR,$  \\
& $h \circ (\id_A, 0\,\circ\,!_A) = f: A \to B,$ $(\init)$ \\
& $h \circ (A \times s) = g \circ h: A \times \N \to B,$ $(\step)$
} 
{ $h = g^\S \circ (f \times \N): 
                             A \times \N \to B \times \N \to B,$
}
in terms of \NAME{Freyd}'s pentagonal diagram:

\begin{minipage} {\textwidth}
$$
\xymatrix@+2em{
& A \times \N 
  \ar [r]^{A \times s}
  \ar [d]_{f \times \N}
  \ar @{} [ddr] |{=}
  \ar @{-->} @/^2pc/ [dd]^h
  & A \times \N 
    \ar [d]_{f \times \N}
    \ar @{-->} @/^2pc/ [dd]^h 
\\
A
\ar [ur]^{(\id,0\,!)}
\ar [dr]_f
\ar @{} [r] |{=}
& B \times \N
  \ar [d]_{g^\S} 
  & B \times \N 
    \ar [d]_{g^\S} 
\\
& B 
  \ar [r]^g
  & B
}
$$
\begin{center} \NAME{Freyd}'s uniqueness \textsc{diagram} $(\FR!)$ 
  \end{center} 
\end{minipage}

\bigskip
Schemata $(\mathrm{it})$ and $(\FR!)$ give in fact the well-known
full schema $(\mathrm{pr})$ of Primitive Recursion:
\inference { (\prr) }
{ $g = g(a): A \to B$ $\mr{PR}$ \emph{(init map)} \\
& $h = h((a,n),b): (A \times \N) \times B \to B$ \emph{(step map)}
}
{ $f = f(a,n) = \prr[g,h]\,(a,n): A \times \N \to B$ in $\PR$ such that \\ 
& $f(a,0) = g(a): A \to B$ \emph{(init)}, and \\
& $f(a,s\ n) = h((a,n),f(a,n)): (A \times \N) \to B,$ $(\step)$ \\
& \qquad as well as \\
& $(\prr!):$ $f$ is \emph{unique} with these properties.
}

This \emph{full schema of $\mr{PR}$} has as consequence in particular
the elegant and powerfull \textbf{Uniqueness schemata} of \NAME{Goodstein,}
$\mr{U_1}$ - $\mr{U_4},$ whith \emph{passive} parameter, 
usually $a \in \N$ here made explicit.  

Using these schemata, \NAME{Goodstein} proves central equations
for \emph{addition, truncated subtraction,} and \emph{multiplication} 
on the ``NNO'' $\N:$

\smallskip
As usual, the basic structure of a \emph{unitary commutative semiring}
on NNO $\N$ -- (plus \emph{truncated subtraction and exponentiation}) -- is
\textbf{defined} -- and \textbf{characterised,} map theoretically -- with
interpretation of Free Variables $a$ and $n$ as \emph{projections}
$a :\,= \ell: A \times \N \to A$ and 
$n :\,= r: A \times \N \to \N$ -- as follows:
\begin{align*}
& 0,\ 1 \bydefeq s\,0: \one \to \N:\ \text{\emph{zero} \resp \emph{one,}} \\
& \text{\emph{addition}}\ a+n = s^\S(a,n): \N^2 \to \N\ \text{by} \\
& \qquad a+0 = a: \N \to \N,\ a+s\,n = s(a+n): \N^2 \to \N, \\
& \text{\emph{truncated subtraction}}\ a \dotminus n: \N^2 \to \N
                                                        \ \text{by:} \\
& \qquad a \dotminus 0 = a,\ a \dotminus s\,n = \pre\,(a \dotminus n), \\
& \qquad \text{with \emph{predecessor}}\ \pre: \N \to \N \\ 
& \qquad \text{PR defined (and characterised) by}
                     \ \pre\,0 = 0,\ \pre\,s\,n = n, \\
& \qquad \text{so for example}
            \ 5 \dotminus 2 = 3,\ \text{but}\ 2 \dotminus 5 = 0; \\  
& \text{\emph{multiplication}}\ a \bs{\cdot} n: \N^2 \to \N\ \text{by} \\
& \qquad a \bs{\cdot} 0 = 0,\ a \bs{\cdot} s\,n = a \bs{\cdot} n + a, \\
& \text{and \emph{exponentiation}}\ a^n: \N^2 \to \N\ \text{by} \\
& \qquad a^0 = 0, \ a^{s\,n} = a^n \bs{\cdot} a.
\end{align*}

\smallskip
\textbf{Remark} on use of free variables:
\ul{Inter}p\ul{retation} of \ul{Free Variables} as (identities) \resp
\emph{projections,} (possibly nested) is subject to formal rules
extending this interpretation in the above example. Vice versa,
projections can be \ul{seen} as free variables, by (re)naming them
with names $a,b,\ldots,x,y,\ldots$ usually standing for free
``individual'' variables. From now on, we will make extensive use 
of Free Variables, \emph{anchored} in the Cartesian structure of 
\emph{basic} theory $\PRa = \PR+(\abstr)$ and its strengthenings.

\smallskip
\textbf{Continuation} of \emph{elementary map equations} for 
NNO $\N:$ \emph{Multiplication} on $\N$ \emph{distributes} 
not only over addition but over \emph{truncated subtraction} 
as well, almost (again) by definition.

\smallskip 
(Boolean) \emph{Logic} and \emph{Order} then are defined, and 
characterised as follows, the latter essentially via truncated 
subtraction: 
\begin{align*}
& \sign = \sign(n): \N \to \N\ \text{by}\ \sign(0) = 0,\ 
                                   \sign(s\,n) = 1 = s\,0: \N \to 2, \\
& \neg = \neg(n): \N \to \N\ \text{by}\ \neg(0) = 1,\ \neg(s\,n) = 0.
\end{align*}
A map $\chi: A \to \N,$ here such a map within $\PR,$ is called a 
\emph{predicate} (on Object $A$) if
  $$\sign \circ \chi = \neg \circ \neg \circ \chi = chi: \N \to \N.$$
By obvious reason, we write for this
  $$\chi: A \to 2 \subset \N.$$
For the moment, this is just notation. Later, in theory $\PRa$ to 
come, Object $2 = \set{n \in \N\,|\,n < s\,s\,0} \subset \N$ will 
become an \emph{Object} on its own right.\footnote{
  The idea to introduce Boolean (!) Free-Varibles predicate 
  Calculus into the Theory of Primitive Recursion this way, 
  without an explicit basic (``undefined'') Object 
  $\two \identic \set{\false,\true} \identic \set{0,1}$
  goes back to \NAME{Goodstein} 1971 and \NAME{Reiter} 1982.}

Using this notation, we get the Boolean Operations, and basic 
binary predicates on $\N$ as
\begin{align*}
& \land = [\,a\,\land\,b\,] \defeq \sign(a \bs{\cdot} b): \N^2 \to \N, \\
& [\,m \leq n\,] \defeq \neg\,[\,m \dotminus n\,]: \N^2 \to 2. \\
& [\,m < n\,] \defeq \sign\,[\,n \dotminus m\,]: \N^2 \to 2, \\
& \qquad \text{and directly from ``\,$\leq$\,'':} \\ 
& [\,m \doteq n\,] = [\,m \leq n\,]\ \land\ [\,n \leq m\,]: \N^2 \to 2.
\end{align*}
The latter predicate is \emph{predicative equality.} It can be extended 
to all \emph{fundamental Objects} (finite bracketed powers of $\N$),
by componentwise conjunctive definition. And furthermore to 
(formal, virtual) predicate \emph{extensions} 
$\set{A\,|\,\chi: A \to 2}$ of $\PR$-predicates $\chi,$ see below.  

\smallskip
$[\,$Addition as well as multiplication -- the latter with arguments 
greater zero -- are strictly \emph{monotonous} with respect to the 
(linear) order introduced above via truncated subtraction.$]$

\smallskip
This \emph{equality ``of individuals''} $[\,a \doteq a'\,]: A^2 \to 2$
is \emph{reflexive, symmetric, and transitive.} 
On $\N$ it satisfies -- with regard to strict order 
$<\,: \N^2 \to \N$ -- the \emph{law of trichotomy.} 

Furthermore it is -- as we expect it from a ``sound'' notion of 
\emph{equality} -- \emph{substitutive} (in \NAME{Leibniz}' sense), 
with respect to the ``basic meta-operations'' \emph{composition,} 
\emph{Cartesian Product,} as well as \emph{Iteration,}
and therefore in particular with respect to \emph{all} operations
(above), introduced on $\N$ via these basic \ul{meta-o}p\ul{erations}. 
 
\smallskip
The key relationship between this \emph{predicate}
and the -- already given -- notion of \ul{e}q\ul{ualit}y 
\ul{between} maps:
\begin{align*}
& f =^{\PR} g: A \to B,\ \text{also written} \\
& \PR \derives\ f = g: A \to B: \\
& \text{Theory}\ \PR\ 
     \text{\ul{derives} \ul{e}q\ul{ualit}y}\ f = g,
\end{align*}
is the following \emph{Equality Definability Schema,} 
which is \ul{derivable} in theory $\PR,$ and also in
Theory $\PRa = \PR+(\abstr)$ to come, as well as in  
all strengthenings of these:

\smallskip
\textbf{Equality Definability Theorem:}
An Arithmetical Theory $\T,$ \ie $\T$ an extension of $\PR,$ admits the 
following schema (EquDef):
\inference { (\mr{EquDef}) } 
{ $f,g: A \to \N$ in $\T,$ \\
& $\T \derives\ \doteq \circ\, (f,g) = true_{\N^2}: \N^2 \to 2$ 
}
{ $\T \derives\ f = g: A \to \N,$ \emph{algebraically:} \\ 
& $f =^\T g: A \to \N.$ 
}
\emph{Equality Definability} extends to $\T$-map pairs 
$f,g: A \to B$ with (common) codomain a (finite) Cartesian Product $B$ 
of Objects $\N,$ a \emph{fundamental} Object -- in $\PR$ -- or even 
$B$ a predicative extension $B = \set{C\,|\,\gamma}$ of such an Object,
an Object of Theory $\PRa$ for short, see below.

\smallskip
\textbf{Remark:} For \textbf{proof} of the laws of multiplication, 
and also for proof of \emph{logically} all important 
\emph{Equality Definability} above, \NAME{Goodstein} 
proves commutativity of the maximum function, namely
  $$\max(a,b) \defeq a+(b \dotminus a) = b+(a \dotminus b) 
                                 \bydefeq \max(b,a): \N^2 \to \N.$$  

We now ``realise'' straight forward the schema of 
\emph{predicate abstraction:} We start with a ($\PR$) predicate 
$\chi: A \to 2$ on a \emph{fundamental} Object $A:$ Object $A$ 
is a (binary bracketed) Cartesian product out of $\N$ and $\one.$ 
Such a predicate $\chi: A \to 2,$ formally: $\chi: \N \to \N$
(see above), can be turned into a \emph{virtual,}
``new'' Object $\set{A\,|\,\chi} = \set{a \in A\,|\,\chi(a)}$ easily:
just take as Objects of \emph{new frame,} extended Theory $\PRa$ the 
fundamental maps of form $\chi: A \to 2$ (\emph{predicates}) within 
\emph{fundamental} Theory $\PR$ above, and as maps between such Objects 
$\set{A\,|\,\chi}$ and $\set{B\,|\,\psi}$ those $\PR$-maps
$f: A \to B$ which transform $\chi$ into $\psi.$ Two such maps
$f,g: A \to B$ are \emph{identified,} declared \emph{equal,} 
\ul{if} they agree ``on'' $\chi.$ This \emph{definitional, conservative}
extension is called here $\PRa = \PR+(\abstr) \bs{\sqsup} \PR,$
the \emph{basic Theory of Primitive Recursion.}


\smallskip
\textbf{Structure Theorem for $\PRa:$} \emph{Basic} Theory $\PRa$ 
becomes a \emph{Cartesian Category} (Theory), has all 
\emph{Extensions} of its predicates and therefore all \emph{equalisers,} 
and contains \emph{fundamental} Theory $\PR$ \ul{embedded} via
  $$\bfan{f: A \to B} \bs{\mapsto} 
       \bfan{f: \set{A\,|\,\true_A} \to \set{B\,|\,\true_B}}.$$
Furthermore, it has all equality predicates -- by restriction --
and admits the schema of \emph{Equality Definability.}

\textbf{Proof} (\NAME{Reiter} 1980): A preliminary version of $\PRa$ 
is constructed as canonical extension of $\PR$'s Class of Objects 
into its predicates and then admitting as \emph{maps} between these 
``new'' Objects those of $\PR$ which are compatible with the 
given, \emph{defining} predicates for these ``Objects''.
This way one gets a Cartesian $\mr{PR}$ theory. Equality of $\PRa$
maps, given by their equality on the predicate of their $\PRa$
Domain Object, is compatible with this Cartesian $\mr{PR}$ structure 
and thereby gives the Cartesian $\mr{PR}$ structure of $\PRa,$ 
defined as Quotient theory by this \emph{notion of equality.} 
Forming the \emph{iterated} is likewise compatible with this
canonical notion of equality in the extended world, and therefore
$\PRa$ becomes this way a $\mr{PR}$ category. By construction,
$\PRa$ gets in addition the wanted predicate-extension and 
equaliser-by-extension structure: it has finitie limits, in
particular pullbacks (and multiple pullbacks)\quad\textbf{\qed}

\bigskip
\textbf{Remarks:}
 
(i) A $\PRa$-map $f: \set{A\,|\,\chi} \to \set{B\,|\,\psi}$
can be viewed as a \emph{defined partial} $\PR$ map from $A$ to $B$
with values in $\psi:$ Set of \emph{defined arguments,} namely
$\set{a \in A\,|\,\chi(a)}$ is PR \emph{decidable.} By \textbf{definitin}
of $\PRa$'s equality, $\PR$-map $f: A \to B$ ``doesn't care'' about
arguments $a$ in the \emph{complement} $\set{a \in A\,|\,\neg\,\chi(a)}.$

So wouldn't it be easier to realise this view to 
\emph{defined partial maps} just by throwing the 
\emph{undefined arguments} into a \emph{waste basket} $\set{\bot}$ say?

But where to place this waste basket, this for each Codomain Object $B?$
The fundamental Objects have a zero-vector as a candidate. For example
we could interprete truncated subtraction as a \emph{defined partial} map
  $$a \dotminus b: \set{(m,n) \in \N \times \N\,|\,m \geq n} \to \N,$$
and throw the complement $\set{(m,n) \in \N \times \N\,|\,m < n}$ into
waste basket $\set{0} \subset \N.$ But this is not a good interpretation
of \emph{truncated} (!) subtraction: Value $0$ is \emph{not} waste, it
has an important meaning as zero.

``The'' waste basket $\set{\bot}$ should be an entity with a
\emph{natural} extra representation, and we should have only 
one such entity in a later theory of defined partial PR maps 
to come. This theory, to be called $\PRaX,$ will be constructed 
with the help of \emph{Universal Object} $\X$ which is to contain 
\emph{codes} of all singletons and (nested) pairs of natural 
numbers, and ``below'' these codes it has room for code of 
\emph{undefined value} symbol $\bott,$ in a ``Hilbert's hotel''.
All this will be carried through within present theory $\PRa.$

\smallskip
(ii) A $\PR$-map $f: A \to B$ such that $f$ ``is'' a $\PRa$-map
\begin{align*}
& f: \set{A\,|\,\chi\,\lor\,\chi': A \to 2} \to \set{B\,|\,\psi},
                             \ \text{also ``works'' as a $\PRa$-map} \\
& f: \set{A\,|\,\chi} \to \set{B\,|\,\psi},\ \text{and a $\PRa$-map} \\ 
& g: \set{A\,|\,\chi} \to \set{B\,|\,\psi\,\land\,\psi'} 
                               \ \text{also ``works'' as a $\PRa$-map} \\
& g: \set{A\,|\,\chi} \to \set{B\,|\,\psi}.
\end{align*}
Since map-properties of \emph{injectivity,} \emph{epi-property}
of $\PR$-maps viewed as $\PRa$-maps, \textbf{depend} on choice 
of hosting (predicative) $\PRa$ Objects -- \textbf{examples} above 
-- \emph{specification} of a $\PRa$ map 
$f: \set{A\,|\,\chi} \to \set{B\,|\,\psi}$ must contain,
besides $\PR$-map $f: A \to B$, Domain and Codomain \emph{Objects} 
$\chi: A \to 2$ and $\psi: B \to 2$ as well. 
This way the members of \ul{Ma}p \ul{set} $\ulfamily$ $\PRa(A,B):\ A,B$ 
$\PRa$-Objects, become mutually \ul{dis}j\ul{oint}. Inclusions
$i: A' \ovs{\subseteq} A''$ are realised in $\PRa$ as restricted
$\PR$-identities 
$\id_A: \set{A\,|\,\chi'} \ovs{\subseteq} \set{A\,|\,\chi''},
                                               \ \chi' \implies \chi''.$

\section{Goodstein Arithmetic $\GA$}

In ``Development of Mathematical Logic'' (Logos Press 1971) 
R. L. Goodstein gives four basic uniqueness--rules for 
Free--Variable Arithmetics. We show here 
that these four rules are sufficient for proving the commutative and 
associative laws for multiplication and the distributive law, for addition 
as well as for truncated subtraction.

\smallskip
We include\footnote{Sandra Andrasek and the author} 
into Goodstein's uniqueness--rules a ``passive parameter'' $a$. 
These extended  rules are provable using Freyd's uniqueness Theorem
$(\mathrm{pr!}),$ part of \emph{full schema} $(\mathrm{pr})$
of Primitive Recursion which he deduces from his uniqueness 
$(\FR!)$ of the \emph{initialised iterated.} \NAME{Freyd} deduces 
the latter from availability of a Natrual Numbers Object $\N$
in \NAME{Lawvere's} sense, and (!) \emph{axiomatic} availability 
of ``higher order'' \emph{internal} $\mathrm{hom}$ objects 
with, again axiomatic, \emph{evaluation} map family for these 
objects, of form $\epsilon_{A,B}: B^A \times A \to B$ within (!) 
the category considered.

\smallskip
\textbf{Goodstein's rules with passive parameter:}
Let $f, g: A \times N \to N$ be primitive recursive maps, 
$s: N \to N$ the successor map $n \mapsto n+1$ and 
$\pre: N \to N$ the predecessor map, usually written as 
$n \mapsto n \dotminus 1.$ 

Then Goodstein's rules read:

\bigskip
$
\begin{array} {l}
\inference{ \mathrm{U_1} }
  { $f(a,sn) = f(a,n): A \times \N \to B$ }
  { $f(a,n) = f(a,0): A \times \N \to B$ } 
\\  \\ 
\bigskip
\inference{ \mathrm{U_2} }
  { $f(a,s\,n) = s\,f(a,n): A \times \N \to \N$ }
  { $f(a,n) = f(a,0)+n: A \times \N \to \N$ }
\\

\bigskip
\inference{ \mathrm{U_3} }
  { $f(a,sn) = \pre\ f(a,n): A \times \N \to \N$ }
  { $f(a,n) = f(a,0) \dotminus n: A \times \N \to \N$ }
\\

\bigskip
\inference{ \mathrm{U_4} }
  { $f(a,0) = g(a,0): A \to \N$ \\ 
  & $f(a,sn) = g(a,sn): A \times \N \to \N$ 
  }
  { $f(a,n) = g(a,n): A \times \N \to \N.$ }
\end{array} 
$

\textbf{Comment:} Theories $\PR$ and $\PRa$ allow, within rules 
$\mathrm{U}_1$ and $\mathrm{U}_4$ above, for replacing $\N$ as a 
Codomain Object, by an arbitrary object $B$ of $\PR$ \resp $\PRa.$ 

Rule $\mathrm{U_4},$ of \emph{uniqueness} of maps defined by
\textbf{case distinction,} is nothing else than \emph{uniqueness} 
of the \emph{induced map out of the sum} 
$A \times \N \iso (A \times \one)+(A \times \N),$ 
this sum canonically realised via \emph{injections} 
$\iota = (\id_A \times 0): A \times \one \to A \times \N$
as well as, right injection: 
$\kappa = \id_A \times s: A \times \N \to A \times \N:$

This uniqueness combined with \NAME{Leibniz}' \emph{compatibility} 
of \emph{induced-map-out-of-a-sum} with map (term) equality, 
compatibility available in Theories $\PR,$ $\PRa,$ and their 
strengthenings. 
 
\smallskip
\textbf{Proof} of these four rules is straight forward for theories
$\PR,\ \PRa$ (and strengthenings), using \NAME{Freyd}'s uniqueness
$(\FR!)$ and uniqueness clause $(\mathrm{pr}!)$ of the 
\emph{full schema of Primitive Recursion} respectively.

\smallskip
$\mr{U_1}$-$\mr{U_4}$ give, by means of a \ul{derived} schema $V_4,$ the wanted

\smallskip
\textbf{Structure Theorem} for NNO $\N:$ $\N$ admits the 
structure of a \emph{commutative semiring} with zero,
\emph{truncated subtraction} $\dotminus\ =\ m \dotminus n \to \N,$
over which multiplication \emph{distributes,} a 
\emph{linear order} $m < n: \N^2 \to 2,$ and \emph{equality predicate}
$\doteq\ =\ [\,m \doteq n\,]: \N^2 \to 2,$ both of the latter 
\textbf{defined} via truncated subtraction. Order and equality 
(predicates) satisfy the \emph{law of trichotomy;} addition is 
strictly monotoneous in both arguments; truncated subtraction is 
weakly monotoneous in first, and weakly antitoneous in second argument, 
wheras multiplication is strongly monotoneous in both arguments, on
$\N_{>0}^2 = \set{\N\,|\,\,> 0}^2.$

\smallskip
\noindent
$\N$ and $\set{n \in \N\,|\,n < 2}$ admit (2-valued) 
\textbf{Boolean Logic} 
$\sign,\ \neg,\ \land,\ \lor,\ \implies.$

\smallskip
\noindent
Last but not least, the maximum $\max(a,b): \N^2 \to \N$ \emph{commutes:}
  $$\PR \derives \max(a,b) \defeq a+(b \dotminus a) 
                   = b+(a \dotminus b) \bydefeq \max(b,a): \N^2 \to \N.$$
This goes into \textbf{Proof} of \textbf{Equality Definability,}
for $\PR,\ \PRa,$ and strengthenings. 
 
\textbf{Proof,} by $\mr{U_1}$ - $\mr{U_4},$ straightforward but tedious.

\section{Theories of Partial $\mr{PR}$ Maps}

We now turn to Extension of $\mr{PR}$ theories -- here Theory $\PRa$
or a strengthening -- into Theories of \emph{partial} $\mr{PR}$ maps 
(\emph{map terms}).
This extension is understood best, wenn looking at the following
\textsc{diagram} which shows \emph{composition} of two such
partial maps, $g: B \parto C$ with $f: A \parto B:$

\begin{minipage} {\textwidth}
$$
\xymatrix@+2em{   
D_h 
\ar @/_2pc/ [dd]_{d_h}^{=} 
\ar [d]^{\pi_\ell} 
\ar [rd]^{\pi_r} 
\ar @/^2.8pc/ [rrdd]^{\widehat{h}} 
  & & 
\\
D_f 
\ar [d]^{d_f} 
\ar [rd]^{\widehat{f}} 
\ar @{} [r] | {\pb} 
& D_g 
  \ar [d]^{d_g} 
  \ar [rd]^{\widehat{g}} 
  \ar @{}[r] | {=} 
  & 
\\
A 
\ar @{^>} [r]^f 
\ar @/_2pc/ @{^>} [rr]_{h\,=\,g\,\widehat{\circ} f}^{\pareq} 
& B 
  \ar @{^>} [r]^g 
  & C 
}
$$
\begin{center} Composition \textsc{diagram} for $\hatS$ \end{center}
\end{minipage}

\bigskip
A partial map \quad
  $f = \bfan{(d_f: \widehat{f}): D_f \to A \times B}: A \parto B$ \\
consists of a $(\mr{PR})$ enumeration $d_f = d_f(\hat{a}): D_f \to A$ 
of \emph{defined-arguments} for $f,$ and a \emph{rule} 
$\widehat{f} = \widehat{f} (\hat{a}): D_f \to B,$ fixing the values
$f(a) \defeq \widehat{f} (\hat{a})$ for \emph{defined} arguments $a \in A,$
\ie $a$ of form $a = d_f(\hat{a}),\ \hat{a} \in D_f.$

Up to here, $f$ defines just a \emph{relation,} in the sense of
\NAME{Brinkmann}/\NAME{Puppe} 1969. What is lacking is
\emph{right uniqueness,} see the following

\bigskip
\textbf{Partial-Map Schema:} 
\inference{ (\hatS) }
{ $\gamma{f} = \gamma{f} (\hat{a}): D_f \to A \times B$\ \ $\bfS$-map, \\
& \qquad called \emph{graph} (of $f: A \parto B$ to be introduced), \\
& $d_f = d_f(\hat{a}) \defeq \ell \circ \gamma{f}: D_f \to A$ 
                           \ \emph{defined arguments enumeration} \\
& $\widehat{f} = \widehat{f} (\hat{a}): D_f \to B$\ \emph{rule} \\
& $\bfS \derives\ d_f(\hat{a}) \doteq d_f(\hat{a}') 
               \implies \widehat{f} (\hat{a}) \doteq \widehat{f} (\hat{a}'):
                 D_f^2 \to 2\ (\hat{a},\ \hat{a}' \in D_f\ free):$ \\
& \qquad \emph{right uniqueness} 
}
{ $f \defeq \bfan{\gamma{f} = (d_f,\widehat{f}): 
                          D_f \to A \times B}: A \parto B$ \\
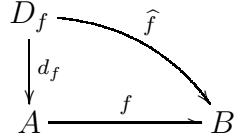
& $\hatS$-morphism, \emph{partial} $\bfS$ map. 
}

\smallskip
\textsc{Diagram} above then shows -- \NAME{Brinkmann}\,
\&\,\NAME{Puppe} type -- \emph{partial map} \emph{composition} 
$g \parcirc f: A \parto B \parto C$ via pullback.

\smallskip
\textbf{Equality} $f \pareq g: A \parto B$ of partial maps 
-- over $\bfS$ -- is established by availability of a pair 
\xymatrix{D_f \ar @<0.4ex> [r]^i & D_{f'} \ar @<0.2ex> [l]^j} 
of \emph{defined-arguments comparison} maps (in $\bfS$) which are
compatible as such with $d_f,d_{f'}$ as well as with 
$\widehat{f},\ \widehat{f'}.$ 
Availability of just on of these $\bfS$ maps, of 
$i: D_f \to D_{f'}$ say, \textbf{defines} ``graph''
\emph{inclusion,} here $f\,\parinc\, f': A \parto B.$   

Basic \emph{compatibility} of (partial) \emph{composition} 
``$\parcirc$'', with graph \emph{inclusion} $\parinc$ -- and 
hence with \emph{partial equality} ``$\pareq$'' then is given 
by the universal properties of (composed) pullback in the following 
\textsc{diagram\footnote{F.\ Hermann}:}      

\bigskip
\begin{minipage} {\textwidth}
$$
\xymatrix{ 
& & & D_{h}
      \ar@{-->} '[lll] '[llldddd]^k [dddd]
      \ar[dr]|{r}
      \ar[dl]|{\ell}
      \ar@/_1.5pc/[ddll]|{d_{h}}
      \ar@/_-1.5pc/[ddrr]|{\widehat{h}}
      \ar@{}[dd]|{\pb} 
\\
& & D_{f}
    \ar@{..>}[dd]^(.3) {i}
    \ar[dr]|{\widehat{f}}
    \ar[dl]|{d_{f}} &
    & D_{g'}
      \ar@{..>}[dd]^(.3)j
      \ar[dr]|{\widehat{g}}
      \ar[dl]|{d_{g}} 
\\
& A
  \ar@{-^{>}}[rr]^(.7){f}_(.7){f'}
  \ar@{-^{>}}@/_1.5pc/[rrrr]^h_{h'} |!{[dr];[rr]}\hole
                                        |!{[drrr];[rr]}\hole   
  & & B
      \ar@{-^{>}}[rr]^(.3)g_(.3){g'}      
      & & C 
\\
& & D_{f}
    \ar[lu]|{d_{f'}}
    \ar[ru]|{\widehat{f'}} 
    & & D_{g'}
        \ar[lu]|{d_{g'}}
        \ar[ru]|{\widehat{g'}} 
        & 
\\
& & & D_{h'}
      \ar[ur]|{r'}
      \ar[ul]|{\ell'}
      \ar@/_-1.5pc/[uull]|{d_{h'}}
      \ar@/_1.5pc/[uurr]|{\widehat{h'}}
      \ar@{}[uu]|{\pb}
}
$$
\begin{center} Compatibility \textsc{diagram} 
                 of $\parcirc$ with $\subseteq$ 
\end{center} 
\end{minipage}
  
\bigskip
Furthermore, composition via pullback above then is associative, by 
\emph{natural equivalence} of the (finite) limits defining
compositions 
  $$h \parcirc (g \parcirc f),\  (h \parcirc g) \parcirc f:
                               A \parto B \parto C \parto D.$$  
\emph{Cylindrification} is componentwise, and gives the Cartesian Product
for $\PRa$ as a \emph{monoidal} -- again bifunctorial one -- within 
extended Theory $\hatPRa \bs{\sqsupset} \PRa.$ But this extended Product
does not have anymore (\NAME{Godement}'s) universal properties of a 
\emph{Cartesian} Product, within $\hatPRa.$

\smallskip 
\emph{Iteration} in $\hatPRa$ works analogeously to composition, using
in this case \emph{pullback iteration.} 

\smallskip
\textbf{Equality Definability:} There is such a Theorem also for 
\emph{partial map} Theory $\hatPRa$ and its strengthenings. 
  
\bigskip
\textbf{Structure Theorem for $\hatS:$}

\begin{enumerate} [(i)]

\item 
$\hatS$ carries a -- canonical -- structure of a 
\emph{diagonal} \emph{symmetric} \emph{monoidal category}, 
with \emph{partial} composition $\parcirc$ and identities introduced above, 
(monoidal) product $\times$ extending $\times$ of ${\bfS},$ \emph{association}
$\ass: (A \times B) \times C \overset{\iso} {\lto} A \times (B \times C),$
\emph{symmetry (``transposition'')} 
$\Theta: A \times B \overset{\iso} {\lto} B \times A,$ and 
\emph{diagonal} $\Delta: A \to A \times A$ inherited from $\bfS;$  
\cf \NAME{Budach} \& \NAME{Hoehnke} 1975 and 
-- later\footnote{there is an earlier preprint of 
\NAME{Budach} \& \NAME{Hoehnke}} -- \NAME{Pfender} 1974 for an axiomatic 
approach to categories with a given type of \emph{substitution}
transformations. Our present theory $\hatS,$ a theory of \emph{partial} 
$\mr{PR}$ maps, is a monoidal category, which has -- in addition to 
\emph{natural} transformations $\ass,\ \Theta,$ and $\Delta$ above --,
so-called \emph{half-terminal} maps, and the former projections as
\emph{half-projective} ones, in the terminology of \NAME{Budach} 
\& \NAME{Hoehnke}, \emph{``half''}  since the latter \emph{natural} 
families of Theory $\bfS,$ are no longer \emph{natural transformations} 
for Theory $\hatS.$ All of this \emph{substitutive} structure is
obviously preserved by the embedding $\bfS \bs{\sqsub} \hatS.$ 

\item 
The defining diagram for an $\hatS$-map -- namely

\begin{minipage} {7 cm}
$$
\xymatrix{ 
D_f
\ar @/^1em/ [rrd]^{\widehat{f}}
\ar [d]^{d_f}
  & 
\\
A 
\ar @{^>} [rr]^f
  & & B
}
$$
\begin{center} Partial Map \textsc{diagram} \end{center} 
\end{minipage}

\bigskip
-- constitutes in fact a \emph{commuting} $\hatS$ diagram. 

Conversely -- with same notation as above -- \textbf{define} the 
\emph{minimised opposite} to $d_f,$ beginning with formally 
\emph{partial,} 
$\hatS$ map 
  $$d'^-_f = \bfan{(d_f,[\ ]_{\widehat{f}}): 
                         D_f \to A \times D_f}: A \parto D_f,$$
\emph{opposite (graph)} to given $\bfS$ map $d_f: D_f \to A.$ 
This opposite is made \emph{right-unique} by selecting 
$D_f$-\emph{minimal} $\widehat{f}$ \emph{equivalence representant}
  $$[\ ]_{\widehat{f}} = [\hat{a}]_{\widehat{f}}
      \defeq \min_{D_f} \set{\hat{a}' \leq \hat{a}\,|\,
                          \widehat{f} (\hat{a}') \doteq_B 
                                         \widehat{f} (\hat{a})}:
                                                      D_f \to D_f:$$
\emph{Minimal} with respect to here \emph{``canonical'',} 
\NAME{Cantor}-ordering on $\bfS$ Object $D_f = \set{D\,|\,\delta: D \to 2}.$ 
This order is \emph{inherited} from $D_f$'s \emph{``mother''} 
fundamental Object, $D,$ say.
This object in turn is (well) ordered via canonical \emph{counting}
  $$\cantor_D = \cantor_D(n): \N \overset{\iso} {\lto} D,$$
(see general schema above of \emph{PR dominated minimum}),
and \textbf{get} the commuting $\hatS$-\textsc{diagram}

\bigskip
\begin{minipage} {5 cm}
\xymatrix{ 
D_f
\ar @/^1em/ [rrd]^{\widehat{f}}
\ar @{} [rd] |{\quad \pareq}  
  & 
\\
A 
\ar @{-^{>}} [u]^{d^-_f}
\ar @{-^{>}} [rr]^f
  & & B
}
\end{minipage}
  \begin{minipage} {3 cm}
  put together:
  \end{minipage}       
\begin{minipage} {5.5 cm}
\xymatrix{ 
D_{f}
\ar @/^1em/ [rrd]^{\widehat{f}}
\ar @<-0.5ex> [d]_{d_f}
\ar @{} [rd] |{\quad \pareq}  
  & 
\\
A 
\ar @{-^{>}} @<-0.5ex> [u]_{d^-_f}
\ar @{-^{>}} [rr]^{f}
  & & B
}
\smallskip
\begin{center} Basic Partial Map \textsc{diagram} \end{center} 
\end{minipage}

\item 
The first factor $f: A \parto B$ in a $\hatPRa$-composition
$h = g \parcirc f: A \parto B \parto C,$ when giving an (embedded)
$\PRa$ map $h: A \to C,$ is itself an (embedded) $\PRa$ map:
\emph{first factor of a total map is total.}

\smallskip
So each \emph{section} (``coretraction'') of $\hatPRa$ is a $\PRa$ map, 
in particular a $\hatPRa$ section of a $\PRa$ map is in $\PRa.$

\item 
$\hatS$ clearly inherits from ${\bfS}$\, 
\NAME{Fourman}'s \emph{uniqueness equation:} 
  $$\text{For}\ h: C \parto A \times B\ \text{in}\ \hatS:
      \ h \pareq (h \parcirc \ell, h \parcirc r): C \parto A \times B,$$ 
where for $f: C \parto A \;,\; g: C \parto B,$
  $$(f,g) \defeq (f \times g) \,\parcirc\, \Delta_C: 
                         C \to C \times C \parto A \times B,$$
with \emph{diagonal} $\Delta_C : C \to C \times C$ of $\bfS.$

\smallskip
This equation guarantees \emph{uniqueness} of the ``\emph{induced}''
$(f,g): C \parto A \times B,$ but $(f,g)$ does not satisfy 
(both of) the Cartesian \ul{e}q\ul{uations} 
  $$\ell \parcirc (f,g) \pareq f 
      \quad \ul{\mr{and}} \quad r \parcirc (f,g) \pareq g,$$
except $f$ and $g$ have \emph{equal domains of definition}, \ie if
$i: D_f \to D_g, \, j: D_g \to D_f$ are available such that 
$d_g \circ\ i = d_f: D_f \to A$ as well as $d_f \circ\ j = d_g: D_g \to A.$


\item 
\emph{Iteration} $f^\S: A \times \N \parto A$ of $\hatS$-endo is
available in $\hatS,$ satisfying the characteristic $\hatS$-equations
\begin{align*}
& f^\S \parcirc (\id_A,0) \bydefeq f^\S \parcirc (A \times 0) \circ \Delta_A
                                      \pareq \id_A: A \to A, \ \ul{and} \\
& f^\S \parcirc (A \times s) \pareq f \parcirc f^\S: A \times \N \parto A.
\end{align*}

\item 
Freyd's uniqueness of the \emph{initialised iterated} holds in $\hatS:$
\inference{ (\FR!)_{\hatS} }
{ $f: A \parto B,\ g: B \parto B,\ h: A \times \N \parto B$ 
                                              in $\hatS$ such that \\
& $h \parcirc (\id_A,0) \parto f: A \parto B$ \ \text{\ul{and}} \\
& $h \parcirc (A \times s) \pareq g \parcirc h: A \times \N \parto B$
}
{ $h \pareq g^\S \parcirc (f \times \N): 
      A \times \N \parto B \times \N \parto B.$ }

$[\,$The latter two statements are not so easy to prove: 
PR \emph{construction} of \emph{comparison} maps is needed, for comparing 
the respective enumerations of defined arguments in the postcedent, 
proceeding from the comparison maps given by the antecedents.$]$ 
 
\item 
For extension $\hatS$ of $\bfS$ again, we get -- by the general
\NAME{Freyd}'s argument -- the corresponding
\emph{full schema of primitive recursion,} namely
\inference{ (\mathrm{pr})_{\hatS} }
{ $g: A \parto B$ in $\hatS$ $(\mathrm{initialisation}),$ \\
& $h: (A \times \N) \times B \parto B$ $(\mathrm{step\ map})$
} 
{ $f = \mathrm{pr} [\,g,h\,]: A \times \N \parto B$ is available in $\hatS,$ \\
& characterised (up to equality $\pareq$) in $\hatS$ by \\ 
& $f \parcirc (\id_A,0) \pareq g: A \parto B$ \text{\ul{and}} \\
& $f \parcirc (A \times s) \pareq h \parcirc (\id_{A \times \N},f)$ \\
& $ \bydefeq 
      h \parcirc ((A \times \N) \times f) \parcirc \Delta_{A \times \N}:$ \\
& $A \times \N \to (A \times \N)^2 \parto (A \times \N) \times B \parto B.$
} 


\end{enumerate}


\smallskip
The \textbf{Proof} of this \textbf{Structure Theorem} is long, 
already since we have to show that many assertions; but mainly 
since assertion (vi) needs some auxiliary arguments. 

Nevertheless, all of these assertions look plausible: they
are ``straightforward'' extrapolations from the case of \emph{finite}
partial maps, by means we have at our disposition for the
\emph{potentially infinite,} primitive recursive case as 
\emph{basic} theory.\footnote{full Proof is ready within 
detailed version, as a pdf file}  


\section{Partial-Map Extension as a Closure Operator} 

\smallskip 
\textbf{Closure:} Theory $\hathatS$ of \emph{partial} maps over $\hatS,$
of \emph{partial} partial maps over $\bfS,$ is (category) equivalent 
to Theory $\hatS.$ Theory $\bfS$ is a strengthening of $\PRa$ as always 
here.

Mutatis mutandis, construction of \emph{Partiality Hull} 
$\hatS \bs{\sqsup} \bfS$ above of a \emph{Cartesian} $\mr{PR}$ theory
$\bfS$ can be applied again to \emph{diagonal monoidal}
Theory $\hatS$ ``again''. Because of lack of Cartesianness this
is more involved, and so is verification of the properties of this
\emph{Double Closure} $\hathatS.$ In particular it
is more difficult to define composition: If you want to go
into this detail, look at next \textsc{diagram:}

\smallskip
For \textbf{defining} composition of such $\hathatS$-morphisms,
composition of say $f: A \pparto B$ and $g: B \pparto C,$ consider the
following $\bfS/\hatS/\hathatS$-\textsc{diagram} which displays the
$\hatS/\bfS$ data of $f$ and $g$ to be \emph{composed} into an
$\hathatS$-morphism $g \pparcirc f: A \pparto B \pparto C:$

\bigskip
\begin{minipage} {\textwidth}
$$
\xymatrix{
D_f
\ar @{-^{>}} [d]_{\gamma{f}}
& D_{\gamma{f}}
  \ar[l]_{d_{\gamma{f}}}
  \ar[dl]_{\widehat{\gamma{f}}}
  & D_{g \pparcirc f}
    \ar @{=} [d]
    \ar[l]_{\pi_{\ell}}
    \ar[r]^{\pi_r}
    \ar[dll]
    \ar[drr]
    & D_{\gamma{g}}
      \ar[r]^{d_{\gamma{g}}}
      \ar[dr]^{\widehat{\gamma{g}}}
      & D_g
        \ar @{-^{>}} [d]^{\gamma{g}}
\\
A \times B
  \ar[d]^{\ell}
  \ar[drr]^{r}
  & & {D_{\gamma{f}} \underset{r \circ \widehat{\gamma{f}}\,,
               \,\ell \circ \widehat{\gamma{g}}} {\times} D_{\gamma{g}}}
      \ar @{} [d]| {\pb}
      & & B \times C
            \ar[dll]_{\ell}
            \ar[d]_{r}
\\
A 
\ar @{-_{>}} [rr]_{f}
& & B
    \ar @{-_{>}} [rr]_{g}
    & & C
}
$$
\begin{center} Composition \textsc{diagram} for $\hathatS$ \end{center}
\end{minipage}

\bigskip
Composition $g \pparcirc f: A \pparto C$ then is \textbf{defined}
to have as \emph{graph} $\gamma_{g \pparcirc f}$ the map 
\emph{``induced''} by the left and right \emph{frame} morphisms of the 
\textsc{diagram,} namely:
  $$\gamma_{g \pparcirc f}
      \defeq (\ell \parcirc \gamma{f} 
                 \parcirc d_{\gamma{f}} \circ \pi_{\ell}\,,
            \,r \parcirc \gamma{g} 
                 \parcirc d_{\gamma{g}} \circ \pi_r):
                      D_{g \pparcirc f} \parto A \times C.$$
The next assertion (really) to be \textbf{proved} is \emph{idempotence}
of our \emph{Closure} operator, namely that each $\widehat{\hatS}$
map $f: A \pparto B$ is represented -- with respect to notion
of equality $\ppareq$ of $\widehat{\hatS},$ by a suitable
$\hatS$-map $h: A \parto B.$

For a \textbf{Proof} look at the following \textsc{diagram,} for 
given $\widehat{\hatS}$-map 
  $$f = \bfan{\gamma{f} 
           = (\ell \parcirc \gamma{f},\ r \parcirc \gamma{f})
             =\,: (d_f,\widehat{f}): 
                D_f \parto A \times B}: A \pparto B:$$

\begin{minipage} {\textwidth}
$$
\xymatrix{
& D_{h} 
  \ar @{=} [d]^{\mathrm{def}}
  \ar @/_2pc/ [ddd]_(0.2){h=\gamma{h}}
  &
\\
& D_{\gamma{f}} 
  \ar @/_2pc/ [lddd]_{d_{h} :\,= \,\widehat{d_f}}
  \ar [d]^{d_{\gamma{f}}}
  \ar @/^2pc/ [dd]^{\widehat{\gamma{f}}}
  \ar @/^2pc/ [rddd]^{\widehat{\widehat{f}}\, =\,: \widehat{h}}
  &
\\
& D_f          
  \ar @{-^{>}} @/^/ [u]^{d^-_{\gamma{f}}}
  \ar @{-^{>}} @/_/ [ldd]_{d_f}
  \ar @{-^{>}} [d]^{\gamma{f}}
  \ar @{-^{>}} @/^/ [rdd]^{\widehat{f}} 
  &
\\
& A \times B 
  \ar [ld]_{\ell}
  \ar [rd]^{r}
  &
\\
A 
\ar @<0.5ex> @{-_{>}} [rr]^{f}
\ar @<-0.5ex> @{-^{>}} @/_1pc/ [rr]^{\ppareq}_{h}
  & {}
    & B
}
$$
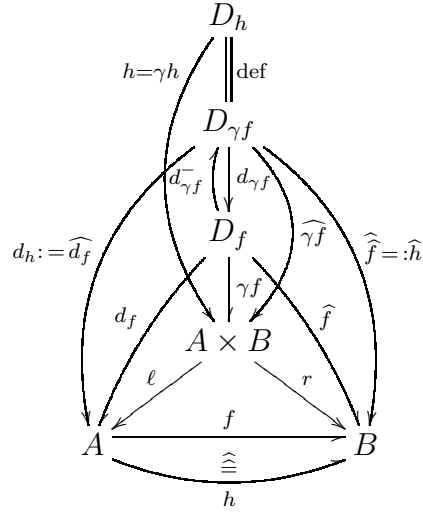
\begin{center} Closure \textsc{diagram} for Extension by partial maps 
  \end{center} 
\end{minipage}

\bigskip
In this \textsc{diagram,} $\gamma{f}: D_f \parto A \times B$ is 
the \emph{graph} of $\hathatS$-morphism $f: A \pparto B$ to be
considered. The $\bfS$-maps $d_{\gamma{f}}: D_{\gamma{f}} \to D_f$ 
(defined-arguments enumeration) and 
$\widehat{\gamma{f}}: D_{\gamma{f}} \to A \times B$ 
(rule) are to define $\gamma{f}: D_f \parto A \times B$
as a \emph{partial} $\bfS$-map, an $\hatS$ morphism.

\smallskip
This \textsc{diagram} shows the way of \textbf{Proof} for

\smallskip
\textbf{Closure Theorem for Extension of Theory $\bfS$ by Partial Maps:}

\emph{Closure} by \emph{Partial Maps} is \emph{idempotent:}
Partial map Closure of theory $\hatS$ is again a diagonal monoidal
category $\hathatS$ which is in fact \emph{equivalent} 
-- as such a category -- to theory $\hatS:$ 
  $$\hathatS \bs{\iso} \hatS.$$
$[\,$Both inherit -- from $\bfS$ -- Object $\N$ as NNO in the sense 
of the (full) schema of $\mr{PR}$ for $\N.]$

\section{$\mu$-Recursion without Quantifiers} 

\textbf{Church type ``Inclusion'':} For given $\PRa$-predicate 
$\ph: A \times \N \to 2,$ \emph{partially defined} ``map''
  $$\mu\ph = \mu\ph(a): A \parto \N,$$ 
classically defined by
$$\mu\ph(a) =
  \begin{cases}
    \min\set{n \in \N\,|\,\ph(a,n)}
                 \ \text{if}\ (\exists n \in \N)\ \ph(a,n) \\
    \text{\emph{undefined}}
                 \ \text{if}\ (\forall n \in \N)\ \neg\,\ph(a,n),
  \end{cases}
$$
has a -- classically \emph{correct} -- \emph{representation} within 
(strengthenings of) Theory $\hatPRa$ as
$$
\xymatrix @+2em{
D_{\mu\ph} \defeq \set{A \times \N\,|\,\ph} = \set{(a,n)\,|\,\ph(a,n)}
\ar[d]_{d_{\mu\ph}\,=\,\ell\,\circ\,\subseteq}
\ar[rd]^{\widehat{\mu}\ph}
&
\\
A
\ar @{-^>} [r]_{\mu\ph}
& \N,
}
$$
Here \emph{defined-arguments} $(\mr{PR})$ enumeration is 
  $$d_{\mu\ph} = d_{\mu\ph} (a,n) \defeq a: 
       \set{A \times \N\,|\,\ph} \subseteq A \times \N 
                                      \overset{\ell} {\lto} A,$$
and rule $\widehat{\mu}\ph: \set{(a,n)\,|\,\ph(a,n)} \to \N$ is 
(totally) $\mr{PR}$ defined by
  $$\widehat{\mu}\ph 
      = \widehat{\mu}\ph(a,n) = \min\set{m \leq n\,|\,\ph(a,m)}: 
                                   \set{(a,n)\,|\,\ph(a,n)} \to \N.$$

\smallskip
\textbf{$\mu$-Lemma:} $\hatS$ admits the 
following (Free-Variables) schema $(\mu),$ operator $\mu$'s 
``property'', combined with \emph{uniqueness} schema $(\mu!),$ 
as a \emph{characterisation} of the $\mu$-operator  
  $\bfan{\ph: A \times \N \to 2} \bs{\mapsto} \bfan{\mu\ph: A \parto \N}$
above: 


\inference { (\mu) } 
{ $\ph = \ph(a,n): A \times \N \to 2 
                       \ \bfS-\text{map (``predicate'')},$ \\
& \qquad $[\,a \in A \ (\free)$ is the \emph{passive} parameter, \\ 
& \qquad $n \in \N \ \free$ the \emph{recursion parameter}$\,]$
}
{ $\mu\ph = \bfan{(d_{\mu\ph}\,,\,\widehat{\mu}\ph): 
                                         D_{\mu\ph} \to A \times \N}:
                                                             A \parto \N$ \\
& \qquad is an $\hatS$-map such that \\ 
& $\bfS \derives\ 
     \ph(d_{\mu\ph} (\hat{a}),\widehat{\mu}\ph(\hat{a})) 
                                       = \true_{D_{\mu\ph}}:
                                                   D_{\mu\ph} \to 2,$ \\
& \qquad $[\,\hat{a} \in D_{\mu\ph}$ free, so just $a \in A$ 
                            of form $a :\,= d_{\mu\ph} (\hat{a})$ \\
& \qquad counts as -- is \emph{enumerated} as -- ``defined argument'' 
                                                      for $\mu\ph\,]$ \\
& + ``argumentwise'' \textbf{minimality:} \\
& $\bfS \derives\ [\,\ph(d_{\mu\ph} (\hat{a}),n) 
                     \implies \widehat{\mu}\ph(\hat{a}) \leq n\,]: 
                                      D_{\mu\ph} \times \N \to 2$
}
as well as \emph{uniqueness,} by \emph{maximal extension}:
\inference { (\mu!) } 
{ $f = f(a): A \parto \N$ in $\hatS$ such that \\
& $\bfS\ \derives\ \ph(d_f(\hat{a}),\widehat{f}(\hat{a}))
                                       = \true_{D_f}: D_f \to 2,$ \\
& $\bfS \derives\ \ph(d_f(\hat{a}),n) 
               \implies \widehat{f}(\hat{a}) \leq n: D_f \times \N  \to 2$
}
{ $\bfS \derives\ f \ \widehat{\subseteq}\ \mu\ph: A \parto \N$ 
                          \ \text{(inclusion of graphs).}
}


\textbf{Proof} of schema $(\mu)$ is by \textbf{definition} of
$\mu\ph = \mu\ph\,(a): A \parto \N.$ \textbf{Proof idea} for
uniqueness schema $(\mu!)$ is ``displayed'' as the following 
\textsc{diagram:}
\bigskip
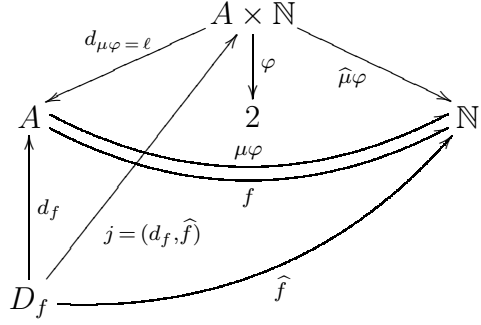
\begin{minipage} {\textwidth}
$$
\xymatrix{
& & A \times \N
    \ar @<-0.5ex> [lld]_{d_{\mu\ph\,=\,\ell}}
    \ar [d]^{\ph}  
    \ar @<0.5ex> [rrd]_{\widehat{\mu}\ph}
    & &
\\
A
\ar @{-^{>}} @<1ex> @/_2pc/ [rrrr]^{\mu\ph}
\ar @{-^{>}} @/_2pc/ [rrrr]_{f}
  & & 2 
      & & \N
\\
& & & &
\\
D_f
\ar [uu]_{d_f}
\ar [rruuu]_(0.3){j\,=\,(d_f,\widehat{f})}
\ar @/_2pc/ [rrrruu]_{\widehat{f}}
  & & & &
}
$$
\begin{center} $\mu$-applied-to-$\bfS$-predicates 
  \textsc{diagram} gives $f\,\parinc\,\mu\ph.$ \end{center}
\end{minipage}

\smallskip
\textbf{Remark:} Within \NAME{Peano}-Arithm\'etique $\PA,$ and hence also
within set theory, our $\mu\ph: A \parto \N$ equals
  $$\mu\ph = \bfan{(\subseteq\,,\widehat{\mu}\ph): 
                 \hat{A} \to A \times \N}: A \supset \hat{A} \to \N,$$ 
with 
 $\hat{A} = \set{\hat{a} \in A\,|\,\exists n \; \ph(\hat{a},n)},$ 
and
 $\widehat{\mu}\ph(\hat{a}) 
   = \min \set{m \in \N\,|\,\ph(\hat{a},m)}: \hat{A} \to \N,$
\ie it is given there by the classical \emph{partial} minimum definition.
But this definition lacks \emph{constructivity,} since 
$\hat{A} \subseteq A$ is in general not $\mr{PR}$ decidable.

\bigskip
\textbf{Conversely,} consider a \emph{partial} $\mr{PR}$ map, 
  $$f = \an{(d_f,\widehat{f}): D_f \to A \times B}: A \parto B,
                                      \ \text{out of}\ \hatPRa.$$
\textbf{Standard, pointed case:} $f$ is \emph{defined} at least 
at one \emph{point,} say at 
$a_0 = d_f(\hat{a}_0): \one \to D_f \to A.$ 

Such $f$ is represented easily, within Theory $\hatPRa,$ by a 
$\mu$-recursive $\hatPRa$-map (followed by a $\PRa$ map), 
namely by
  $$\mu[\!\![f]\!\!] \defeq 
      (\widehat{f} \circ \cnt_{D_f}) \parcirc \mu\ph_f: 
                               A \parto \N \parto D_f \to \N.$$
$\PRa$-predicate $\ph_f = \ph_f (a,n): A \times \N \to 2$ 
is given as
  $$\ph_f = \ph_f(a,n) 
      \defeq [\,a \doteq_A d_f \circ \cnt_{D_f} (n)\,]: 
                     A \times \N \to A \times D_f \to 2\ \mr{PR}.$$
Here $\cnt_{D_f}: \N \to D_f$ designates a (retractive) 
$(\mr{PR})$ \emph{count} of $D_f.$ For disposing on this \emph{count} 
of $D_f,$ we had to assume that $D_f$ comes with a $\PRa$-\emph{point,} 
$\hat{a}_0: \one \to D_f$ above. 

Partial map $\mu\ph_f: A \parto \N$ is the genuine $\mu$-recursive
kernel of $\mu$-representation $\mu[\!\![f]\!\!]: A \parto B$
of (pointed) partial map $f: A \parto B.$ We count 
\emph{composition} of $\mu$-recursive maps with $\mr{PR}$
maps equally under the $\mu$-recursive ones. So in this
sense, $\mu[\!\![f]\!\!] \pareq f: A \parto B$ is a 
$\mu$-recursive representant of $f$ within $\hatPRa$ and
its strengthenings, a (partial) map in $\muR,$ and in 
$\mu\bfS$ for strengthenings $\bfS$ of $\PRa.$ 

\smallskip
The case that $D_f$ has \emph{no point,} and is nevertheless 
\emph{not} $\bfS$-\ul{derivabl}y empty, causes a 
\textbf{problem,} formally.
We ``solve'' this problem by modifying (extending) the 
definition of $\mu[\!\![f]\!\!]: A \parto \N$ as follows:

$\PRa$-Object $D_f$ is predicative restriction $D_f: D \to 2$ of
a \emph{fundamental} Object $D,$ which comes as such with a
(componentwise) zero $0: \one \to D$ as (privilegded) \emph{point.}
This Object admits a \NAME{Cantor} count $\cnt_D: \N \ovs{\iso} D.$
(Trivial exception: for $D = \one,\ \cnt_D =\ !: \N \to \one$ 
is still a retraction.)

In present general case, replace -- on the way -- Object 
$A \subset \X$ by sum, (disjoint) union 
$\one \oplus A \subset \X,$ and define  
  $\rph_f = \rph_f(a,n): A \times \N \to 2$ by the following
$\PRa$-\textsc{diagram:}

\bigskip
\begin{minipage} {\textwidth}
\xymatrix{
A \times \N
\ar[rr]^{\rph_f}
\ar[d]^{\subset}_{\iota \times \N}
\ar @{} [ddrr]| \defeq
& & 2
\\
(\one \oplus A) \times \N
\ar[d]^{(\one \oplus A) \times \cnt_D}
& & (\one \oplus A) \times (\one \oplus A)
    \ar[u]_{\doteq}
\\
(\one \oplus A) \times D
\ar @{=} [rr]
& & (\one \oplus A) \times ((D \smallsetminus D_f) \oplus D_f)
    \ar [u]^{\id}_{\times (!_{D \smallsetminus D_f} \oplus d_f)}
}
\end{minipage}

\bigskip
$\rph_f = \rph_f(a,n): A \times \N \to 2$ in Free-Variables notation:
\begin{align*}
& A \times \N \owns (a,n) 
                \overset{\iota} {\mapsto} (a,n) 
                  \mapsto (a,\cnt_D(n)) \\
& \mapsto 
  \begin{cases}
  \false\ \myif\ \cnt_D(n) \notin D_f,\ \text{(``outside'' case)}, \\
  [\,a \doteq_A d_f(\cnt_D(n))\,]\ [\ \in 2\ ]
      \ \myif\ \cnt_D(n) \in D_f. 
  \end{cases}
\end{align*}  
This given, we \textbf{define} in this general case 
  $$\mu[\!\![f]\!\!] \defeq 
      \widehat{f} \circ \cnt_D \parcirc \mu\rph_f: 
                            A \parto \N \to D \to B.$$
Note first that $\widehat{f}: D_f \to B$ comes by (formal)
Domain restriction of a genuine $\PR$ map 
$\widehat{f}: D \to B':$ $D_f = \set{D\,|\,D_f},$ 
$D$ (and $B'$) fundamental: This by definition of 
\emph{maps} of Theory $\PRa = \PR+(\abstr).$

Second: wider count $\cnt_D,$ available
in particular for $D$ as a fundamental Object,
Codomain-restricts here nicely, gives ``again'' 
$\mu$-representation of $f,$ here via 
$D \supseteq D_f = \set{D\,|\,D_f: D \to 2}.$ 

\smallskip
This taken together gives $\mu$-representation of general 
partial $\mr{PR}$-map \\ 
$f: A \parto B$ as 
  $$f \pareq \mu[\!\![f]\!\!] 
        \bydefeq (\widehat{f} \circ \cnt_D) \parcirc \mu\rph_f: 
    \xymatrix{A \ar @{-^>} [r]^{\mu\rph_f} 
              & \N \ar[r]^{\cnt_D}
                & D \ar[r]^{\widehat{f}} & B.}$$ 
\smallskip
So we have reached

\smallskip
\textbf{Another Proved Instance of Church's Thesis:} 
\begin{itemize}
\item
The notion of a $\mu$-\emph{recursive} (partial) map is 
\emph{equivalent} to that of a \emph{Partial} $\mr{PR}$ map, 
over ``all'' Theories of Primitive Recursion. 

\item
Theories $\hatPRa$ and $\muR$ are equivalent, and the 
\textbf{Closure Theorem} $\widehat{\hatS} \bs{\iso} \hatS$
above then shows $\mu\muR \bs{\iso} \muR:$ 

\item
\emph{Level-one $\mu$-recursion is enough} for getting all 
$\mu$-recursive maps. By the above, this gives the well known 
corresponding result for $\while$ p\ul{ro}g\ul{rams}: 
\emph{one $\while$ loop is enough:} Any such 
program can be equivalently transformed into a $\while$ loop 
program without nesting of $\while$ loops.

\item
All this works as well for \emph{strengthenings} $\bfS,\ \hatS$
of $\PRa$ and $\hatPRa$ respectively. We would name the 
corresponding Theory of $\mu$-recursion $\mu\bfS \bs{\iso} \hatS.$
\end{itemize} 

\bigskip
\textbf{Conclusion} so far:
\begin{itemize}

\item 
We can \emph{eliminate formal existential quantification}
 -- as well as (individual, formal) \emph{variables} -- 
from the theory of $\mu$-recursion: we interprete application
of $\mu$-operator to predicates of theories $\bfS$ strengthening
$\mr{PR}$ Theory $\PRa = \PR+(\abstr)$ as suitable \emph{partial}
maps, maps in Theory $\hatS.$

\item 
The $\mu$-operator canonically extends to all \emph{partial}
predicates $\ph: A \times \N \parto 2,$ and associates 
to them just partial maps $\mu\ph: A \parto \N,$ within $\hatS$ 
itself. So, ``once again'', we see, that theories $\hatS$ of 
\emph{partial $\mr{PR}$ maps} are \emph{closed,} this time under the 
$\mu$-operator, ``in parallel'' to \emph{Closure} 
of $\hatS$ under forming \emph{partial} maps: 
\emph{partial partial $\mr{PR}$ maps} ``are'' \emph{partial} $\mr{PR}$ maps.

\item
We have the following chain of categorical \emph{equivalences}
of theories considered so far:
  $$\bfS \bs{\sqsubset} \mu\bfS \bs{\iso} \mu\mu\bfS
      \bs{\iso} \hathatS \bs{\iso} \hatS \bs{\sqsupset} \bfS,$$
the inclusions being \emph{diagonal-monoidal $\mr{PR}$ compatible} with
the equivalences.

$[\,$A \emph{partial} $\mr{PR}$ map $f: A \parto B$ which is, 
``by hazard'', a \emph{total} map -- discussion of overall 
\emph{termination} = \emph{total definedness} in part RCF 2 
($\eps\&\mathcal{C}$), is in general \emph{not} itself $\mr{PR}$: 
only its graph $(d_f,f): D_f \to A \times B$ is $\mr{PR}$. 
\NAME{Ackermann} type maps, in particular \emph{evaluation} of 
all $\mr{PR}$-map-codes, are formally partial maps. In 
\emph{well defined cases,} they can be \emph{forced} -- by 
plausible additional \textbf{axiom} -- to become 
\emph{on-terminating,} \ie to get defined-argument enumeration
epimorphic.$]$   

\item Conversely, application of the $\mu$-operator, already just to
$\PRa$-predicates $\ph = \ph(a,n): A \times \N \to 2,$  
\emph{generates} all \emph{partial} $\PRa$-maps $f: A \parto B.$ 


\item As important special cases of basic $\mr{PR}$ theories $\bfS$ 
we have at the moment Theory $\PRa = \PR+(\abstr)$ itself as well as
the $\mr{PR}$ \emph{trace} $\PA \restr \mathrm{PR_A}$ of $\PA:$
All $\PRa$-maps (map terms) with all those equations in between, 
which are derivable by $\PA:$ Our theories, notions, and results have a 
structure-preserving \ul{Inter}p\ul{retation} into (within) 
Peano-Arithmetic $\PA.$

Same for \textbf{set theories} in place of $\PA,$ in particular
$\ZF,\ \ZFC$ and their first order subsystems 
$\oneZF,\ \oneZFC = \oneZF+\mathbf{ACC}.$ \ $\mathbf{ACC}$ is the 
Axiom of \emph{Countable Choice.}
\end{itemize}

\section{Content Driven Loops, 
                                  in particular $\bf{\while}$ Loops} 

Classically, \emph{with} variables, a \ul{while} loop 
$\wh = \wh\,[\,\chi\,|\,f\,]$ is ``defined'' in \emph{pseudocode} by
\begin{align*}
  \wh(a) :\,= [\,& a' :\,= a; \\
                 & \while\ \chi(a')
                     \ \ul{\mathrm{do}}\ a' :\,= f(a') \ \ul{\mathrm{od}}; \\
                 & \wh(a) :\,= a'\,].
\end{align*}
The formal version of this -- within a \emph{classical,} element based
setting -- is the following partial-(\NAME{Peano})-map characterisation:
$$
\wh(a) = \wh\,[\,\chi\,|\,f\,]\,(a) = 
\begin{cases}
  a\ \text{if}\ \neg\,\chi(a) \\
  \wh(f(a))\ \text{if}\ \chi(a)
\end{cases} 
: A \parto A.
$$
It is possible to give a \emph{static} \textbf{Definition} of 
$\wh = \wh\,[\,\chi\,|\,f\,]: A \parto A,$ within 
$\hatPRa\ \bs{\sqsupset}\ \PRa$ (and strengthenings) as follows:  
\begin{align*}
& \text{With}\ \ph = \ph_{\,[\,\chi\,|\,f\,]} (a,n) 
                                \defeq \neg\,\chi\ f^n(a) \\
& \qquad\qquad \bydefeq \neg\,\chi\,f^\S(a,n): 
                            A \times \N \to A \to 2 \to 2, \\
& \text{the}\ \while\ \text{loop} \\ 
& \qquad \wh = \wh\,[\,\chi\,|\,f\,]: A \parto A\ \text{is given as} \\
& \wh \defeq f^{\mu\,\ph_{\,[\,\chi\,|f\,\,]}} (a) \\ 
& \bydefeq f^\S \parcirc (\id_A,\mu\,\ph_{\,[\,\chi\,|f\,\,]}) \\
& \bydefeq f^\S \parcirc (A \times \mu\,\ph_{\,[\,\chi\,|\,f\,]}) 
                                              \parcirc \Delta_A: \\
& A \to A \times A \parto A \times \N \to A. 
\end{align*}
In this generalised categorical sense, we have within theories 
$\hatS\ \bs{\sqsupset}\ \bfS$ ($\bfS$ strengthening $\PRa$), all 
$\while$ loops.


\bigskip
\textbf{Characterisation Theorem} for $\while$ loops \emph{over} $\bfS,$
within Theory $\hatS:$
For $\chi: A \to 2$ (\emph{control}) and $f: A \to A$ (\emph{step}),
both -- for the time being -- $\bfS$-maps,
\ $\while$ loop $\wh = \wh\,[\,\chi\,|\,f\,]: A \parto A$ 
(as defined above), is characterised by the following 
\emph{implications} within $\hatS:$
\begin{align*}
& \hatS\ \derives\ 
    \neg\,\chi \circ a \implies \wh \parcirc a \doteq a: A \parto 2, 
                                                         \ \text{ul{and}} \\
& \hatS\ \derives\ 
    \chi \circ a \implies \wh \parcirc a \doteq \wh \parcirc f \circ a. 
\end{align*}


\textbf{\emph{Dominated} Termination} of $\while$ loop
  $$\wh = \wh[\chi\,|\,f] = \wh[\chi\,|\,f]\,(a): A \parto A,$$
at argument $a \in A,$ is expressed by
  $$[\,m\ \deff\ \wh[\chi\,|\,f]\,(a)\,] \defeq \neg\,\chi\ f^m(a):
                                                   \N \times A \to 2:$$
In words: \emph{iteration of endo $f: A \to A,$ applied to argument $a,$ 
reaches \emph{stop} condition $\neg\,\chi$ after (at most) $m$ 
steps.} Here both: argument $a \in A$ as well as iteration counter
$m \in \N,$ are free variables (categorically: projections).  
$[\,m\ \deff\ \wh[\chi\,|\,f]\,(a) \doteq_A \bar{a}\,]
                                        \ \text{is to mean:}$
\begin{align*}
& [\,m\ \deff\ \wh[\chi\,|\,f]\,(a)\,] \\ 
&   \ \land\ \wh[\chi\,|\,f]\,(a) 
              \doteq_A f^{\min\set{n \leq m\,|\,\neg\,\chi\,f^n(a)}} (a) \\
& \bydefeq f^\S(a,\min\set{n \leq m\,|n\,\neg\,\chi\,f^n(a)}) 
                                                      \doteq_A \bar{a}: \\
& \N \times A^2 \to 2,\  
                 m \in \N,\ a,\bar{a} \in A\ \free.
\end{align*}
$[$``Variable'' $n \in \N$ used in the $\min$ $\mr{PR}$ Operator is \emph{auxiliary,}
\emph{bounded} by $m.$ It does not count as a (free) variable.$]$

From a \emph{logical} point of view, there are -- at least -- the
following open \textbf{Questions,} in

\bigskip
\textbf{Arithmetics Complexity Problem:} 
\begin{enumerate} [(i)]
\item
Does Theory $\PR$ admit \emph{strict, consistent}
strengthenings, or is it a \emph{simple theory,} will say that it
admits its given notion of equality and the indiscrete (inconsistency)
equality of its maps as only ``congruences?'', \cf a simple \emph{group} 
which has as \emph{normal subgroups} only itself and $\set{1}.$
Because of reasons to be explained in later work, my guess here is:
$\PR$ \emph{admits} non-trivial strengthenings, in particular
I suppose that the $\mr{PR}$ \emph{trace} of $\PA,$ explained above, 
is a strict strengthening of $\PR$ \resp $\PRA = \PR+(\abstr).$ 

We cannot exclude at present that all these strengthenings of $\PR$
make up a whole \ul{lattice} of (Free-Variables) Arithmetical Theories, 
each of them giving particular, ``new'' features to Primitive Recursive Arithmetics.

\item
Already at start we possibly have such a strengthening: If 
Free-Variables (``Free Variables'' in the classical sense) 
\emph{Primitive Recursive Arithmetic} $\PRA$ is \textbf{defined}
to have as its terms all map terms obtainable by the (full) schema of 
Primitive Recursion, and as formulae just the \emph{defining equations} 
for the maps introduced by that schema, then I see no way to 
\textbf{prove} all of the usual semiring equations for $\N:$
 
We \emph{need} Freyd's \textbf{uniqueness} $(\FR!)$ -- section 1 above -- of
the \emph{initialised iterated:} From this \NAME{Horn} clause we can
show (!) in particular \NAME{Goodstein}'s uniqueness rules
$U_1$ to $U_4$ on which \emph{his} Proof of the semiring properties of $\N$
is based. He takes these rules as \textbf{axioms.}

My guess is here -- if I have understood right the definition 
of $\PRA$ -- that $\PR = \PRA+(\FR!)$ is a strict strengthening
of $\PRA,$ at least if there is no ``underground'' connection to
the set theoretic view of maps as (possibly infinite)
\emph{argument-value tables.}

\item
So again, Arithmetic would \emph{become} simpler, if Theory $\PR$ 
would turn out to be \emph{simple.} If not, we have a diversity 
of Arithmetic\textbf{s,} a diversity intuitively
far below such issues as Independence of the Axiom of Choice or
of the Continuum Hypothesis. At least the latter is open in the context
of a \textbf{Constructive Analysis} based on map theoretic, Free-Variables
(variable-free) Primitive Recursive and $\mu$-recursive Arithmetics.
\end{enumerate}   

  


  \newpage

\section*{References}



  
\smallskip
  
  \quad\NAME{J.\ Barwise} ed. 1977: \emph{Handbook of Mathematical Logic.}
  North Holland.
  
  \NAME{H.-B.\ Brinkmann, D.\ Puppe} 1969: 
  \emph{Abelsche und exakte Kategorien, Korrespondenzen.} 
  L.N. in Math. \textbf{96.} Springer.

  \NAME{L.\ Budach, H.-J.\ Hoehnke} 1975: \emph{Automaten und Funktoren.}
  Akademie-Verlag Berlin.


 
  \NAME{C.\ Ehresmann} 1965: \emph{Cat\'egories et Structures.} Dunod Paris.

  \NAME{H.\ Ehrig, W.\ K\"uhnel, M.\ Pfender} 1975: Diagram Characterization
  of Recursion. LN in Comp. Sc. 25, 137-143.

  \NAME{S.\ Eilenberg, C.\ C.\ Elgot} 1970: \emph{Recursiveness.}
  Academic Press.
  
  
  \NAME{S.\ Eilenberg, S.\ Mac Lane} 1945: General Theory of Natural 
  Equivalences. \emph{Trans.\  AMS} 58, 231-294.
  
  
  \NAME{P.\ J.\ Freyd} 1972: Aspects of Topoi. 
  \emph{Bull.\ Australian Math.\ Soc.} \textbf{7,} 1-76.
  
  

  \NAME{R.\ L.\ Goodstein} 1971: \emph{Development of Mathematical
  Logic,} ch. 7: Free-Variable Arithmetics.  Logos Press.
  
  \NAME{F.\ Hausdorff} 1908: 
  Grundz\"uge einer Theorie der geordneten Mengen. 
  \emph{Math.\ Ann.} \textbf{65}, 435-505. 



  \NAME{D.\ Hilbert}: Mathematische Probleme. Vortrag Paris 1900. 
  \emph{Gesammelte Abhandlungen.} 
  Springer 
  1970. 
  
  
  \NAME{P.\ T.\ Johnstone} 1977: \emph{Topos Theory.} Academic Press

  \NAME{A.\ Joyal} 1973: Arithmetical Universes. Talk at Oberwolfach.
  
  
  \NAME{J.\ Lambek, P.\ J.\ Scott} 1986: \emph{Introduction to higher order 
  categorical logic.} Cambridge University Press.
  
  
  \NAME{F.\ W.\ Lawvere} 1964: An Elementary Theory of the Category of
  Sets. \emph{Proc.\ Nat.\ Acad.\ Sc.\ USA} \textbf{51,} 1506-1510.

  \NAME{F.\ W.\ Lawvere, S.\ H.\ Schanuel} 1997: 
  \emph{Conceptual Mathematics, A first introduction to categories.}
  Cambridge University Press.



  \NAME{S.\ Mac Lane} 1972: \emph{Categories for the working mathematician}. 
  Springer.
  
  
  \NAME{B.\ Pareigis} 1969: \emph{Kategorien und Funktoren}. Teubner.
  
  \NAME{R.\ P\'eter} 1967: \emph{Recursive Functions}. Academic Press.

  \NAME{M.\ Pfender} 1974: Universal Algebra in S-Monoidal Categories.
  Algebra-Berichte Nr. 20, Mathematisches Institut der Universit\"at
  M\"unchen. Verlag Uni-Druck M\"unchen.


  \NAME{M.\ Pfender} 2008 RCF1d: Theories of PR Maps and Partial PR Maps,
  detailed version. pdf file.





  \NAME{M.\ Pfender, M.\ Kr\"oplin, D.\ Pape} 1994: Primitive
  Recursion, Equality, and a Universal Set. 
  \emph{Math.\ Struct.\ in Comp.\ Sc.\ } \textbf{4,} 295-313.
  
  

  \NAME{H.\ Reichel} 1987: \emph{Initial Computability, Algebraic
  Specifications, and Partial Algebras.} Oxford Science Publications.

  \NAME{W.\ Rautenberg} 1995/2006: \emph{A Concise Introduction to 
  Mathematical Logic.} Universitext Springer 2006.

  \NAME{R.\ Reiter} 1980: Mengentheoretische Konstruktionen in arithmetischen
  Universen. Diploma Thesis. Techn.\ Univ.\ Berlin.

  \NAME{R.\  Reiter} 1982: Ein algebraisch-konstruktiver Abbildungskalk\"ul 
  zur Fundierung der elementaren Arithmetik. Dissertation, rejected by the 
  Mathematics Department of TU Berlin.
  
  \NAME{L.\  Rom\`an} 1989: Cartesian categories with natural numbers object.
  \emph{J.\  Pure and Appl.\  Alg.} \textbf{58,} 267-278.
  



  \NAME{W.\ W.\ Tait} 1996: Frege versus Cantor and Dedekind: on the concept
  of number. Frege, Russell, Wittgenstein: \emph{Essays in Early Analytic
  Philosophy (in honor of Leonhard Linsky)} (ed. W.\ W.\ Tait). Lasalle:
  Open Court Press (1996): 213-248. Reprinted in \emph{Frege: Importance 
  and Legacy} (ed. M.\ Schirn). Berlin: Walter de Gruyter (1996): 70-113.

  \NAME{U.\ Thiel} 1982: Der K\"orper der rationalen Zahlen im
  arithmetischen Universum. Diploma Thesis. Techn.\ Univ.\ Berlin.

  \NAME{A.\ Tarski, S.\ Givant} 1987: \emph{A formalization of set theory 
  without variables}. AMS Coll.\ Publ.\ vol.\ 41.

\bigskip

\bigskip



  \noindent Address of the author: \\
  \NAME{M. Pfender}                       \hfill D-10623 Berlin \\
  Institut f\"ur Mathematik                              \\
  Technische Universit\"at Berlin         \hfill pfender@math.TU-Berlin.DE\\
  

\vfill

\end{document}